\documentclass[a4paper,fleqn]{cas-dc}

\usepackage[numbers, sort&compress]{natbib}
\usepackage{geometry}

\usepackage{subcaption}
\usepackage{subcaption}
\usepackage{natbib}
\usepackage{adjustbox}
\usepackage{caption}
\usepackage{hyperref}
\usepackage{bm}
\hypersetup{
    colorlinks=true,
    linkcolor=blue,
    urlcolor=blue,
    citecolor=blue,
}

\begin{document}
\let\WriteBookmarks\relax
\def\floatpagepagefraction{1}
\def\textpagefraction{.001}
\shorttitle{Data Paper}
\shortauthors{Abdelmuhsen et~al.}

\title [mode = title]{Numerical modeling of pavement mechanical behavior
under the TSD} 
\author[1]{A. ABDELMUHSEN}[orcid= 0000-0002-7077-0591]
\cormark[1]
\address[1]{Univ Gustave Eiffel, MAST-LAMES, Nantes campus, F-44344 Bouguenais, France, abdelmuhsen.abdelgader@univ-eiffel.fr}

\address[2]{Univ Gustave Eiffel, MAST-EMGCU, F-77454 Marne-la-Vallée Cedex 2, France, franziska.schmidt@univ-eiffel.fr}
\author[1]{J-M. SIMONIN}
\author[2]{F. SCHMIDT}
\author[1]{D. LIÈVRE}
\author[1]{A. COTHENET}
\author[1]{M. FREITAS}
\author [1] {A. IHAMOUTEN}

\cortext[cor1]{Corresponding author}

\renewcommand{\thefigure}{\centering\arabic{figure}}

\begin{abstract}
The Subgrade Resilient Modulus (\(M_R\)) is a fundamental metric for assessing the performance of pavement foundations, as it reflects the bearing capacity and stiffness of subgrade materials. Traditionally, \(M_R\) has been estimated through backcalculation methods using deflection measurements from Falling Weight Deflectometers (FWD). However, these conventional methods often encounter limitations that can affect both their accuracy and efficiency. To address these issues, we propose a novel approach utilizing an inverse model that substitutes the traditional FWD method with Traffic Speed Deflectometer (TSD) deflection slope (\(D_S\)) measurements for estimating \(M_R\). This new method employs Machine Learning (ML) techniques to improve the precision and computational efficiency of \(M_R\) estimations. Our approach involves constructing a synthetic dataset through a numerical forward model to develop an inverse model capable of estimating \(M_R\) from \(D_S\) data. This data paper presents the creation of the forward model dataset, which includes numerically simulated TSD deflection slope measurements. The forward model is based on Burmister’s elastic linear isotropic theory and is implemented using advanced techniques in Alizé-Lcpc software. The dataset integrates various material and structural parameters to accurately reflect diverse pavement behaviors. Deflection data are generated under dual-wheel loading conditions, with key deflection and slope parameters calculated using the superposition principle and first derivative methods.
\end{abstract}


\maketitle

\section{Introduction}
The Subgrade Resilient Modulus ($M_R$) of road structure is a critical parameter in evaluating pavement foundation performance, as it directly measures subgrade materials' bearing capacity and stiffness. Traditionally, $M_R$ is estimated using backcalculation methods based on Falling Weight Deflectometer (FWD) deflection measurements. However, conventional techniques face operational limitations that may compromise the accuracy and efficiency of these estimates \cite{ref12, ref13}.\\To overcome these challenges, the authors propose an inverse model that replaces the traditional FWD-based approach with the Traffic Speed Deflectometer (TSD) deflection slope ($D_S$) for estimating $M_R$. This approach leverages Machine Learning (ML) models to enhance the precision and computational efficiency of $M_R$ predictions. The proposed model uses a synthetic dataset to build through a numerical forward model, designed to create an inverse model that estimates $M_R$ using $D_S$ data \cite{Abdelmuhsen23a, Abdelmuhsen24, Abdelmuhsen22b, Abdelmuhsen23b}.\\Thus; this data paper outlines the development of the forward model dataset, which consists of 235 numerically simulated TSD deflection slope measurements. The forward model is based on Burmister’s \cite{ref7} elastic linear isotropic theory, and the simulations are performed using advanced techniques in the Alizé-Lcpc software \cite{ref14}. The dataset encompasses various material and structural properties, comprehensively representing various pavement behavior factors. The deflection data is calculated using the TSD under dual-wheel loading conditions for a pavement structure. Key deflection and slope parameters are derived by applying the superposition principle and the first derivative method.

\section*{TSD Principle of Operation}   
\medskip
The TSD uses an articulated truck with a rear axle applying a load of around 100 kN to the pavement while operating at speeds of 70-80 km/h. Over ten Doppler laser sensors measure deflection velocity near the load, as direct measurements at the load point are impractical \cite{ref18, ref3}.\\To estimate the deflection slope, these sensors capture vertical deflection velocity, which is corrected for vehicle speed to obtain deflection slope. Then the slope is converted into pavement deflection through curve fitting or numerical integration \cite{ref18,ref16,ref5,ref4}. The resulting deflection basin data allows the extraction of parameters like the Basic Damage Index (BDI) and Surface Curvature Index (SCI) to assess the pavement's structural health and residual life \cite{ref19}.  

\begin{figure}[h!]
    \centering
    \includegraphics[width=9cm]{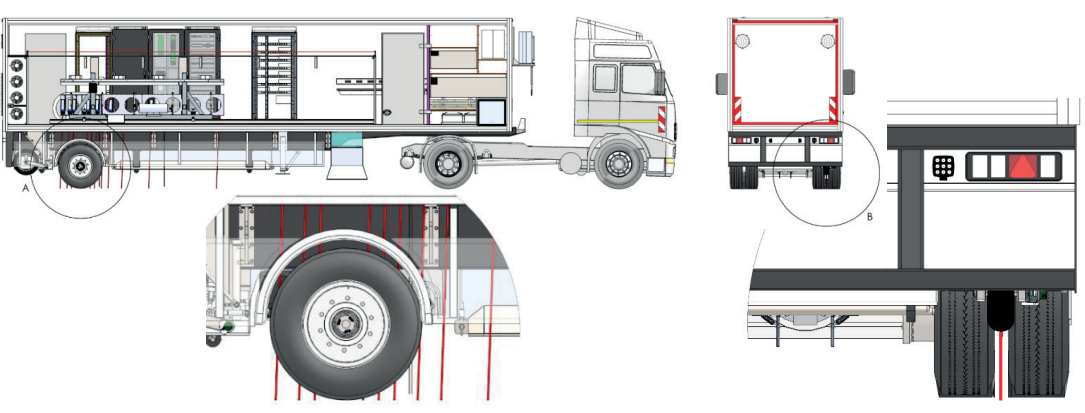}
   
 \caption{TSD measurment system}
    \label{[1]}
\end{figure}

\section*{Numerical Forward Modelling}

\subsection {Alize-LCPC}
Alize-LCPC is the reference software for pavement design in France. It incorporates Burmister’s multi-layer elastic linear model to compute pavement stresses, strains, and deflections. The design process compares calculated values across pavement layers to allowable limits based on material fatigue and rutting behavior, adhering to French standards (NF P98-086/2019). Pavement material properties are sourced from a built-in database \cite{ref12, ref13,ref7,ref14}.

\subsection {Database Parameters}
The input data consists of two main components: Pavement structure and TSD load configuration.

\subsection{Pavement structure: Mechanical properties}
\begin{figure}[h!]
  \centering
\includegraphics[width=0.5\linewidth]{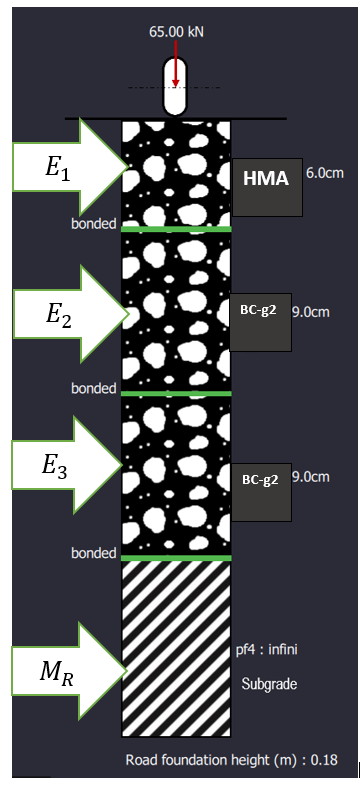}
  \caption{Structure of the pavement under study \cite{ref13}. In scientific terminology from  \cite{ref12} , ($HMA$) denotes a Hot Mix Asphalt. Similarly, ($BC-g2$) represents grade 2 aggregate bituminous concrete, while ($pf_4$) refers to class 4 subgrade.}
  \label{fig:1}
\end{figure}
\begin{table}[h!]
  \centering
  \renewcommand{\arraystretch}{2.5}
  \caption{The characteristics of the pavement structure utilized to derive the forward model (The French catalog \cite{ref12,ref13})}\label{tab:1}
  \begin{adjustbox}{max width=1\linewidth}
  \begin{tabular}{@{}ccccccr@{}}
    \toprule
    Layer (i) & Thickness (m) & Material & $E_i$ (MPa) & $\nu$ & Interface & $^\circ$C \\
    \midrule
    1 & 0.06 & HMA & 7000 & 0.35 & (1) bound & 15 \\
    2 & 0.09 & BC-g2 & 9300 & 0.35 & (1) bonded & 15 \\
    3 & 0.09 & BC-g2 & 9300 & 0.35 & (1) bonded & 15 \\
    \rowcolor{gray!70}
    \textbf{4} & $\bm{\infty}$ & \textbf{Pf\_4} & \textbf{16-250} & \textbf{0.35} & \textbf{(1) bonded} & \textbf{15} \\
    \bottomrule
  \end{tabular}
  \end{adjustbox}
\end{table}
\begin{itemize}
  \item \textbf{Thickness (m):} Represents the thickness of each pavement layer.
  \item \textbf{Material:} The type of material used in each layer, such as Hot Mix Asphalt (HMA) and Bituminous Concrete (BC-g2).
  \item \textbf{Modulus ($E_i$):} Indicates the elastic modulus (in MPa) for each material, reflecting its stiffness.
  \item \textbf{Poisson's Ratio ($\nu$):} A value representing the material's deformation behavior under load.
  \item \textbf{Interface:} Refers to the bonding condition between layers (e.g., bonded or unbonded).
  \item \textbf{Temperature ($^\circ$C):} The temperature at which the material properties are specified.
\end{itemize}
As illustrated in Figure \ref{fig:1} and Table \ref{tab:1}, this study focuses on the modulus of resilience, specifically denoted as $M_{R}$. Accurate estimation of $M_{R}$ is critical for evaluating pavement performance. To achieve this, a robust numerical forward model was developed, which involved simulating data for the variable $S_{V}$ across a range of $M_{R}$ values. This simulated dataset serves as the foundation for training the machine learning (ML) model, which subsequently estimates $M_{R}$ based on Traffic Speed Deflectometer (TSD) measurements.\\The numerical forward model was constructed through an analytical-driven case study that employed predefined parameters (feature engineering) to build a controlled synthetic dataset. This dataset was meticulously designed to facilitate a thorough analysis and accurate categorization of the $M_{R}$ values. Table \ref{tab:1} details the characteristics of this dataset.\\The data, as presented in Figure \ref{fig:1}, demonstrates a variation in the modulus of resilience $M_R$ ranging from 16 MPa to 250 MPa. Initially, the modulus values for all other layers were maintained constant to ensure consistency throughout the study. This approach also allowed for an examination of the sensitivity and impact of each layer’s modulus on the soil independently.\\Furthermore, the temperature across all layers was held constant at 15°C to mitigate the effects of temperature fluctuations on the results.

\begin{table*}[t]
\centering
\caption{Coordinates of TSD Simulated Measurements with Load Values per Wheel}
\label{tab:2cc}
\resizebox{\textwidth}{!}{%
\begin{tabular}{c | *{11}{c}}
\hline
Load       & L1    & L2    & L3    & L4    & L5    & L6    & L7    & L8    & L9    & L10   \\
\hline
$R_x$ (m)  & 0     & 0     & 0     & 0     & 8.150 & 8.150 & 8.150 & 8.150 & 11.750 & 11.750 \\
$R_y$ (m)  & -0.187 & 0.187 & 1.913 & 2.287 & -0.187 & 0.187 & 1.913 & 2.287 & -0.187 & 2.287 \\
$F$ (ton)  & 2.875 & 2.875 & 2.875 & 2.875 & 1.55  & 1.55  & 1.55  & 1.55  & 3.15  & 3.15  \\
\hline
\end{tabular}
}
\end{table*}

\section{TSD Loading mechanism}
The following tables provide an overview of key aspects of the Traffic Speed Deflectometer (TSD) mechanical system. They encompass operational parameters, simulated measurement coordinates, and sensor placements \cite{ref18,ref16,ref5,ref4}.\\Understanding these details is crucial for both establishing the database parameters and accurately interpreting TSD data. This information represents the input components necessary for database utilization and is essential for evaluating pavement performance effectively. Accurate knowledge of these parameters ensures reliable results and facilitates meaningful analysis of pavement behavior and conditions.
\subsection{Coordinates of TSD Simulated Measurements}

Table~\ref{tab:2cc} details the coordinates of the TSD simulated measurements along with the load values per wheel. This table is crucial for understanding the spatial arrangement of the load points and the corresponding sensor measurements.

\begin{itemize}
  \item \textbf{Vertical Distance ($R_x$):} The vertical distance from the reference point of the load to the measurement sensor, recorded in meters.
  \item \textbf{Horizontal Distance ($R_y$):} The horizontal distance from the load point to the sensor, recorded in meters.
  \item \textbf{Load per Wheel ($F$):} The load applied by each wheel, measured in tons.
\end{itemize}

\subsection{TSD Load Configuration}
Table~\ref{tab:4aa} presents the configuration parameters of the TSD load. These parameters define the loading characteristics of the TSD system, which are critical for accurate pavement assessment. 

\begin{table}[h!]
\centering
\renewcommand{\arraystretch}{2}
\caption{TSD Load Configuration}
\label{tab:4aa}
\begin{tabular}{|c|c|c|c|c|c|}
\hline
Parameter & $F$ & $P$  & $V_h$  & $R$  & $f$ \\
\hline
Value     & 32.5     & 0.92      & 20          & 0.15    & 10       \\
\hline
\end{tabular}
\end{table}

\begin{itemize}
  \item \textbf{Axle Load ($F$):} Represents the total load applied by the TSD axle in kilonewtons (kN).
  \item \textbf{Contact Pressure ($P$):} The pressure exerted at the contact point between the tire and the pavement in megapascals (MPa).
  \item \textbf{Speed ($V_h$):} The constant speed of the TSD vehicle in meters per second (m/s).
  \item \textbf{Tire Radius ($R$):} The radius of the TSD tire in meters.
  \item \textbf{Load Frequency ($f$):} The frequency at which the TSD system applies the load, measured in Hertz (Hz).
\end{itemize}

\subsection{Laser Doppler Sensor Position}

Table~\ref{tab:3} provides the positions of the Laser Doppler sensors relative to the TSD load. This information is essential for understanding the spatial arrangement of the sensors and their role in capturing deflection measurements.

\begin{table}[h!]
\centering
\renewcommand{\arraystretch}{2}
\caption{TSD Laser Doppler Sensor Position}
\label{tab:3}
\begin{tabular}{c | c | c | c | c | c | c | c | c}
\hline
Sensor   & $S_1$ & $S_2$ & $S_3$ & $S_4$ & $S_5$ & $S_6$ & $S_7$ & $S_8$ \\
\hline
$S_x$ (m) & 0.1   & 0.2   & 0.30  & 0.45  & 0.60  & 0.90  & 1.1   & 3.8   \\
\hline
\end{tabular}
\end{table}

\begin{itemize}
  \item \textbf{Sensor Position ($S_x$):} The horizontal location of each Laser Doppler sensor relative to the TSD load, measured in meters.
\end{itemize}

The following section will delve into the data processing workflow, outlining the protocol and results associated with each file in the dataset. 

\section{Structured Data Processing Framework}
The synthetic dataset was constructed using Alize-LCPC software, a proprietary tool developed by the University of Gustave Eiffel. For over 30 years, Alize-LCPC has been the leading industrial reference for pavement modeling and design in France, despite not being open-source. This software is protected under the title "Alizé-Lcpc" and version "ALIZE V2 Base 2.0," and is certified by the Agence pour la Protection des Programmes (APP) \cite{ref14} .
\begin{itemize}

    \item \textbf{Alizé-Lcpc Input File:}  
    
    The data processing workflow begins with an input file for Alizé-Lcpc, containing the TSD load configuration and parameters for the target pavement structure. This file includes 235 distinct soil modulus values to simulate deflection slope under varying soil conditions.  
     Below is a screenshot of Alizé-Lcpc input file.
    
    \begin{figure}[h!]
        \centering
        \includegraphics[width=1 \linewidth]{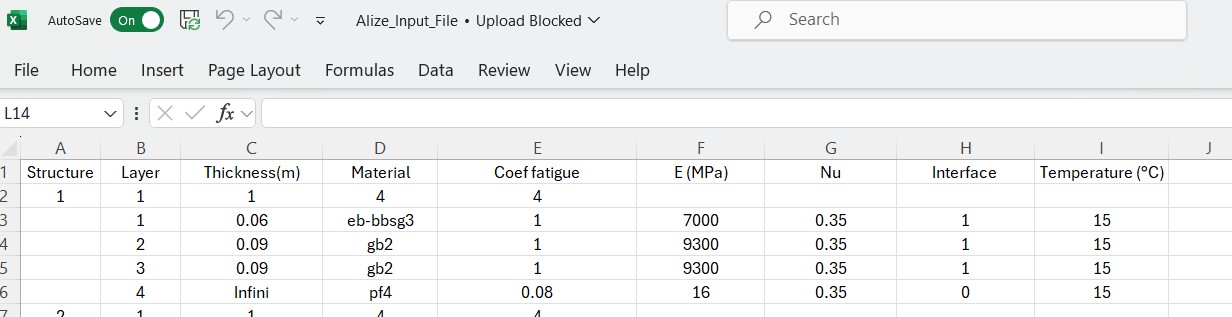} 
        \caption{Screenshot of the Alizé-Lcpc Input File}
        \label{fig:alize_input}
    \end{figure}

    \item \textbf{Alize Output Files:}  
Alizé-LCPC generates 235 output files, each representing pavement deflection under dual-wheel loading for a specific soil modulus from the input. These files are concatenated into a single file, with each column corresponding to a distinct soil modulus value. Figure 3 provides a screenshot of the Alizé-LCPC output files.

    \begin{figure}[h!]
        \centering
        \includegraphics[width=1 \linewidth]{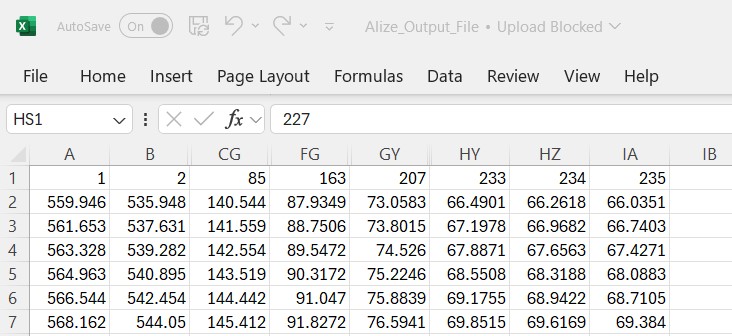} 
        \caption{Screenshot of the Alizé-Lcpc Output Files}
        \label{AO}
    \end{figure}

    \item \textbf{Python: Slope Simulation:}  
    This file contains the simulated deflection slope data from a Python-based simulation program. It includes data for 235 soil modulus values, representing deflections under the TSD. The simulation applies both the superposition principle and the first derivative method. The file comprises columns for the horizontal offset from a 3.8-m reference point, total deflection beneath the TSD's rear wheels (µm), and the horizontal hypothesis distance (m) from the load to each sensor. Additionally, it includes individual wheel deflections (µm), aggregated deflection (using the superposition principle), and the rate of deflection slope change per unit length (µm/m). A corrected deflection slope is provided for increased accuracy at 3.8 m from the load point: figure 4, a screenshot of the Python simulation file.

    \begin{figure}[h!]
        \centering
        \includegraphics[width=1 \linewidth]{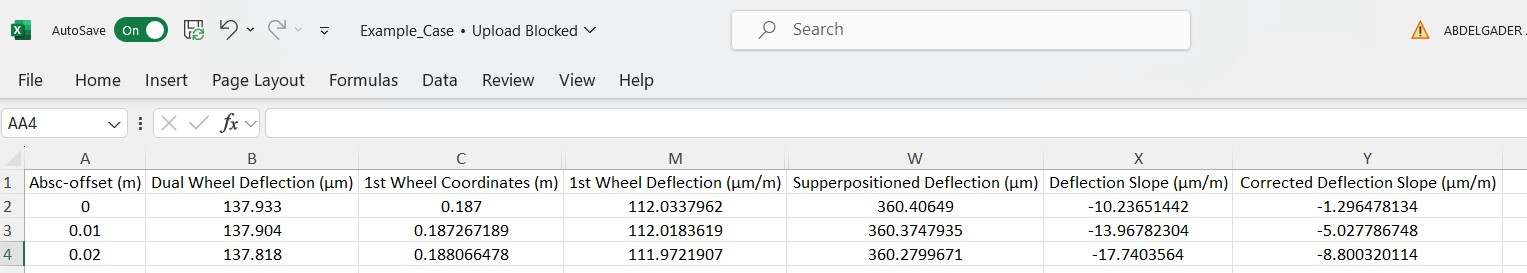} 
        \caption{Screenshot of the Python Slope Simulation File}
        \label{Python}
    \end{figure}

    \item \textbf{Python: Concatenated Deflection Slope Database File:}  
    This final file compiles and consolidates the 235 individual deflection slope files simulated in Python. It represents the Concatenated deflection slope data for all 10 wheels of the TSD, computed using the superposition principle and the first derivative method in the previous step. The figure below shows a screenshot of the concatenated deflection slope database file.

    \begin{figure}[h!]
        \centering
        \includegraphics[width=1 \linewidth]{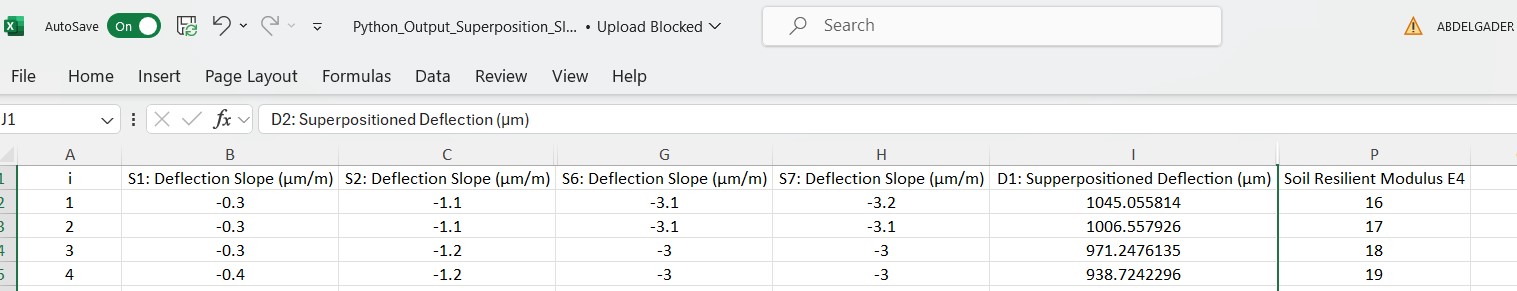} 
        \caption{Screenshot of the Python Concatenated Deflection Slope Database File}
        \label{fig:alize_input}
    \end{figure}

\end{itemize}

To this end, combining these files and the processing steps provides a detailed and integrated view of the relationship between soil stiffness and pavement deflection behavior.

\section{Results}
\subsection{ Mechanical Behaviour Modelling }
\begin{figure*}[h]
\centering
\subfloat[]{\includegraphics[width=5.5 in]{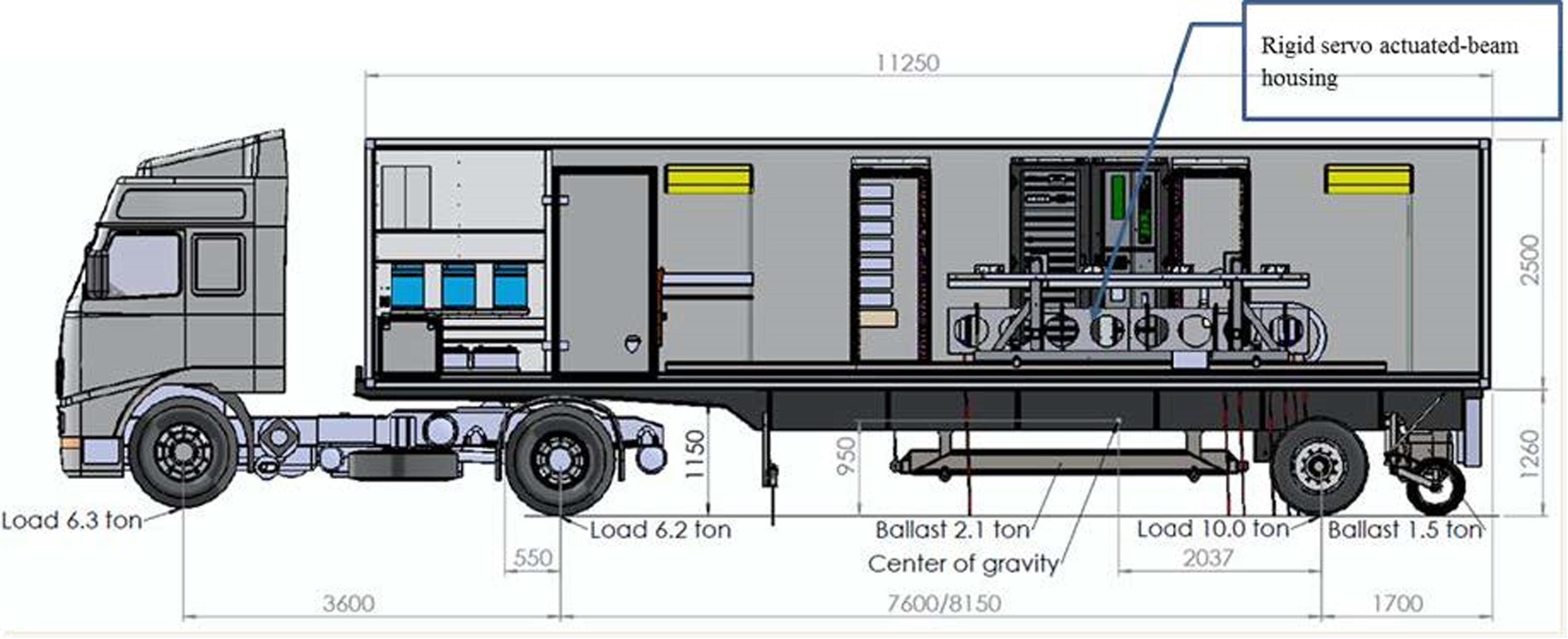}%
\label{subfig:2a}}
\hfil
\subfloat[]{\includegraphics[width=5.8in]{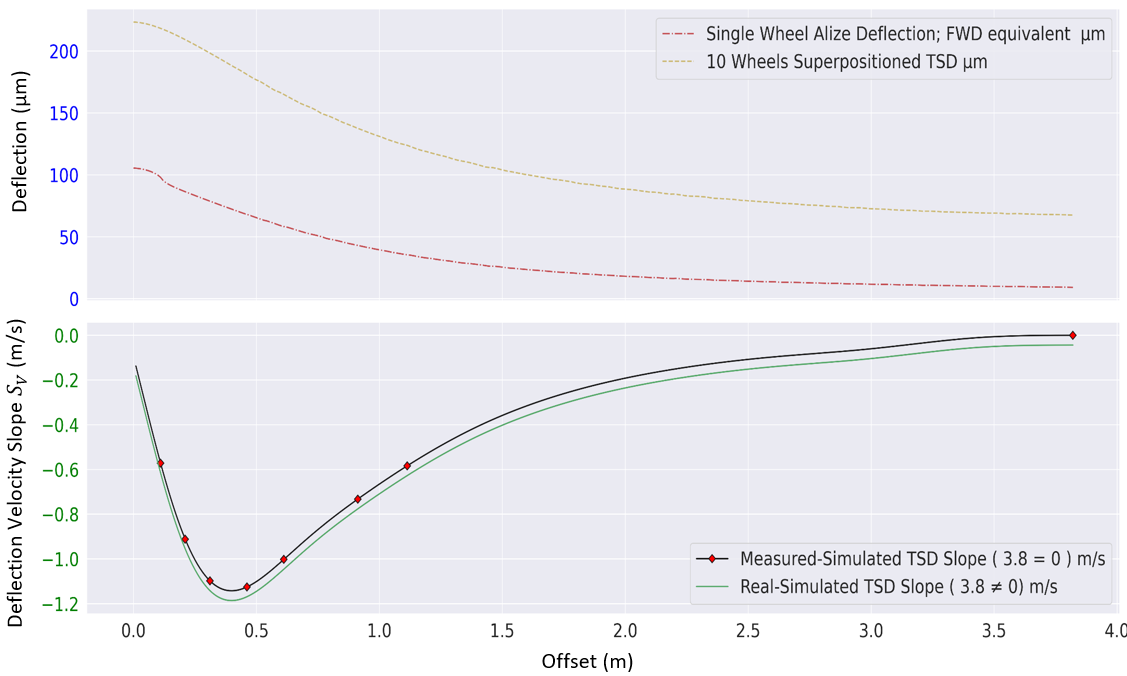}%
\label{subfig:2b}}
\caption{TSD Numerical Model: (a) TSD load configuration \cite{ref16}, (b) Simulation of pavement behavior under TSD  \cite{Abdelmuhsen23a, Abdelmuhsen24, Abdelmuhsen22b, Abdelmuhsen23b}.}
\label{fig:2}
\end{figure*}
To simulate pavement behavior under the TSD system, a load superposition principle was employed. This method is appropriate as it aligns with the deflection behavior hypothesis used in this study, based on Burmister's linear elastic isotropic model. This approach is commonly adopted in pavement engineering to assess soil behavior under load conditions \cite{Abdelmuhsen23a, Abdelmuhsen24, Abdelmuhsen22b, Abdelmuhsen23b}.\\The TSD system is composed of ten wheels ($Fig$ \ref{fig:2}\subref{subfig:2a}), and the simulation involves calculating the total deflection by summing the contributions from all TSD wheels \cite{ref3}. After simulating the overall deflection behavior of the TSD system as illustrated in ($Fig$ \ref{fig:2}\subref{subfig:2b}), the $S_{V}$ parameter is derived by calculating the first-order partial derivative of the combined deflection. This step models the rate of change in pavement deflection under the entire TSD load. Both the deflection and $S_{V}$ curves are presented in ($Fig$ \ref{fig:2}\subref{subfig:2b}).\\To proceed with this, the configuration of the TSD load and sensor positions must be defined, as shown in ($Fig$ \ref{fig:2}\subref{subfig:2a}). Seven laser Doppler sensors ($Sn_{1}$ to $Sn_{7}$) are positioned relative to the TSD wheel loads, and these sensors measure $S_{V}$ beneath the loads.\\To simplify the TSD calibration process, a zero value is artificially assigned to the reference sensor $Sn_{8}$, located 3.8 meters from the load point \cite{ref16}. However, this assumption introduces a systematic bias in the TSD measurements. The impact of $Sn_{8}$ on measurement uncertainty is demonstrated in Figure \ref{fig:2}\subref{subfig:2b}, where the difference between the green curve ($Sn_{8}$ $=$ 0) and the black curves ($Sn_{8}$ $\neq$ 0) illustrates its effect on TSD measurement precision.\\Additionally, ($Fig$ \ref{fig:2}\subref{subfig:2b}) shows that deflection is calculated at 1 cm intervals across every 4 m offset, a practice supported by prior research \cite{ref18}. These aspects are crucial in improving the understanding of the TSD numerical configuration and the factors that influence pavement deflection simulations.

\subsection{Data Modelling and Processing}
\begin{figure}[h!]
\centering
\includegraphics[width=1 \linewidth]{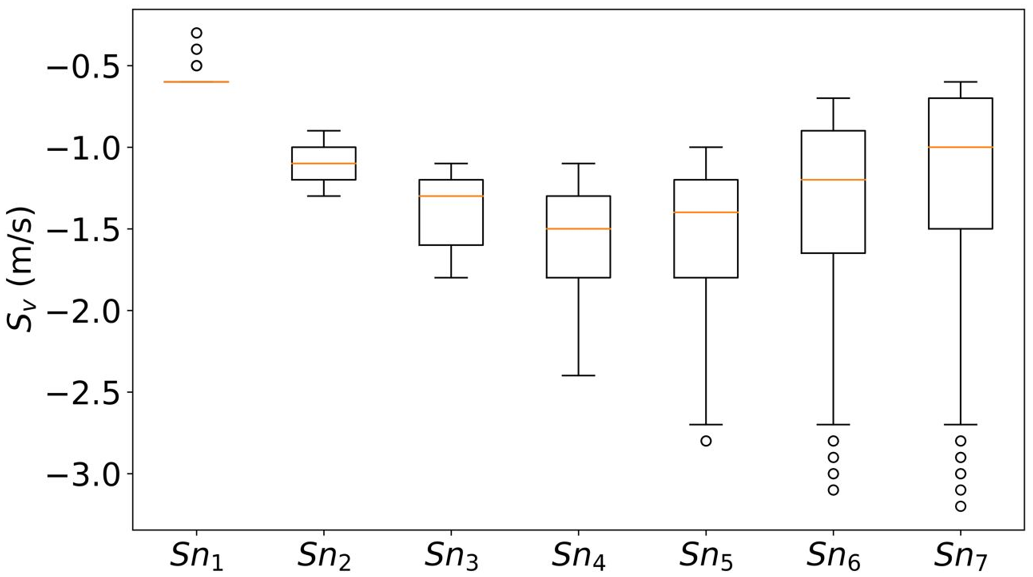}
\caption{$Sn_{1}$ to $Sn_{7}$ Deflection slope analysis \cite{Abdelmuhsen23a, Abdelmuhsen24, Abdelmuhsen22b, Abdelmuhsen23b}}
\label{fig:3}
\end{figure}

This section outlines the process of numerical data collection and describes the key input features used in the study. The dataset comprises 235 instances of simulated $S_{V}$ values, which serve as equivalent representations for varying $M_R$ modulus values ranging from 16 to 250 MPa. The variables $Sn_{1}$ to $Sn_{7}$ correspond to the simulated TSD laser Doppler deflection velocity slopes, with $Sn_{1}$ being the sensor closest to the load application point and $Sn_{7}$ being the farthest. The target variable, $M_R$, refers to the soil's modulus of elasticity and is predicted based on the values recorded by $Sn_{1}$ through $Sn_{7}$.\\Analysis of Figure \ref{fig:3} shows that the variability in $S_{V}$ increases at positions farther from the loading center. Therefore, sensors located at the edges of the deflection basin capture more detailed insights into the soil's mechanical behavior and load-bearing capacity. This observation aligns with the findings of \cite{ref19}, which highlight the correlation between deflection bowl characteristics and mechanistic evaluations of pavement structures. These parameters are instrumental in identifying sections with differing structural capacities and performance levels.

\section{Conclusion}

\begin{itemize}
    \item \textbf{Numerical Forward Model:} A robust numerical forward model was developed using Burmister's elastic linear isotropic theory, implemented via the Alizé-Lcpc software. This model generated simulated deflection slope measurements under dual-wheel loading conditions, providing valuable data for analyzing pavement behavior.
    
    \item \textbf{Dataset Creation:} The study established a detailed dataset incorporating material and structural parameters, which is essential for estimating the Subgrade Resilient Modulus ($M_R$) from Traffic Speed Deflectometer (TSD) measurements. The dataset ranges from 16 MPa to 250 MPa, with consistent parameters for accurate analysis.
    
    \item \textbf{TSD Measurement and Loading Mechanism:} Key aspects of the TSD's mechanical system, including load configuration, sensor positioning, and deflection measurement, were meticulously simulated. This involved defining the TSD axle load, contact pressure, and sensor placements to ensure accurate data collection.
    
    \item \textbf{Sensor Analysis:} Analysis of deflection slopes from multiple laser Doppler sensors revealed that sensors positioned farther from the load point provide more detailed insights into the mechanical response of the pavement. This enhanced variability supports more accurate estimation of the soil modulus.
    
    \item \textbf{Data Processing and Simulation:} A machine learning model can be effectively trained using the synthetic dataset generated from the forward model to estimate $M_R$ values based on deflection slope data. This approach enhances the ability to evaluate pavement performance and structural integrity based on TSD measurements.
\end{itemize}

\section*{Data Availability: \href{https://doi.org/10.57745/QNUCI6}{DOI},  License: \href{https://spdx.org/licenses/etalab-2.0.html}{Etalab 2.0}}
The dataset utilized in this research will be openly accessible to all researchers, accompanied by a comprehensive description in a data paper. It has been assigned an official DOI: \href{https://doi.org/10.57745/QNUCI6}{https://doi.org/10.57745/QNUCI6} and is distributed under the license: \href{https://spdx.org/licenses/etalab-2.0.html}{Etalab 2.0}.

\section*{Conflict of Interest}
The authors declare that they have no financial or personal conflicts of interest that could have influenced the work presented in this paper.

\section*{Acknowledgment}
The first author, on behalf of all co-authors, expresses gratitude to Région Pays de la Loire (France) for their support in co-financing the PhD thesis, which this paper is part of.

\newpage

\bibliographystyle{cas-model2-names}


\end{document}